\batchmode
\documentclass[11pt]{article}
\usepackage{epsfig}
\begin{document}

\title{On the Stochastic Kuramoto-Sivashinsky Equation}
\author{Jinqiao Duan\thanks{Partially supported by the National Science
Foundation Grant DMS-9704345.}
 \ and Vincent J. Ervin \\
Department of Mathematical Sciences \\ 
Clemson University \\ 
Clemson, S.C. 29634-1907 \\ 
U.S.A.  }
\date{ }
\maketitle

\newcommand{\integer}{\mbox{${\rm Z\!\!Z}$}}
\newcommand{\real}{\mbox{${\rm I\!R}$}}
\newcommand{\prob}{\mbox{${\rm I\!P}$}}
\newcommand{\calF}{\mbox{${\mathcal F}$ }}
\newcommand{\calE}{\mbox{${\mathcal E}$ }}
\newcommand{\calT}{\mbox{${\mathcal T}$ }}
\newcommand{\ep}{\mbox{${\epsilon}$}}
\newcommand{\bfx}{\mbox{${\bf x}$}}
\newcommand{\bfu}{\mbox{${\bf u}$}}
\newcommand{\bfv}{\mbox{${\bf v}$}}
\newcommand{\bfn}{\mbox{${\bf n}$}}
\newcommand{\bfb}{\mbox{${\bf b}$}}
\newcommand{\bfr}{\mbox{${\bf r}$}}
\newcommand{\bfxi}{\mbox{${\bf \xi}$}}
\def\qed{\hbox{\vrule width 6pt height 6pt depth 0pt}}
   
\newcommand{\e}{\epsilon}
\renewcommand{\a}{\alpha}
\renewcommand{\b}{\beta}
\renewcommand{\aa}{\mbox{$\cal A$}}
\renewcommand{\L}{\mbox{$\cal L$}}
 
\renewcommand{\theequation}{\thesection.\arabic{equation}}
\newtheorem{corollary}{Corollary}[section]
\newtheorem{lemma}{Lemma}[section]
\newtheorem{theorem}{Theorem}[section]
\catcode`\@=11
\@addtoreset{equation}{section}

\begin{abstract}
In this article we study the solution of the Kuramoto--Sivashinsky
equation on a bounded interval subject to a random forcing term. We
show that a unique solution to the equation exists for all time and
depends continuously on the initial data.

 
\end{abstract}

\section{Introduction}
\label{sec-intro}

In this paper we investigate the existence and uniqueness of 
the solution to the  
Kuramoto-Sivashinsky (K--S) equation subject to a random 
forcing term. Specifically, the solution of
\begin{equation}
du  + (u_{xxxx} +u_{xx} + u u_x) dt - dw =0,
\label{stkeq}
\end{equation}
where  $w$ is a Q--Wiener process in a probability space
$(\Omega, {\cal F}, \prob)$. 
The Wiener process $w$ takes value 
in a Hilbert space to be specified later.
The distributional derivative of $w(t)$ represents an external random force.

The usual K--S equation ((\ref{stkeq}) without the $dw$ term) has 
been studied as a prototypical example for an infinite dimensional
dynamical system. It possesses a finite dimensional maximal
attractor (\cite{nic891,Collet,Ilynko,Goodman}) and 
inertial manifold (\cite{Foias,Jolly,Temam_Wang,Robinson}).

Equation (\ref{stkeq}) arises in the modeling of 
surface erosion via ion sputtering in amorphous materials
\cite{cue951}. The random forcing term in the model accounts for the
fluctuations in the flux of the bombarding particles.

Herein we confine our attention to the case of $u$ restricted to the
interval $I := (-l , l)$, subject to the given initial condition $u_{0}$, and
homogeneous Dirichlet boundary conditions, i.e.
\begin{equation}
   u(0 , x) \ = \ u_{0}( x ) \, , -l < x < l \, , \mbox{ and }
     u(t , -l) \ = \ u(t, l) \ = 0 \, \mbox{ for } t > 0 .
\label{bc1}
\end{equation} 
We show that for any $T > 0$ there exists a unique solution
 to (\ref{stkeq}),(\ref{bc1}) for  $0< t < T$, 
and establish a priori estimates for the solution. The 
approach we follow is similar to that for establishing existence and uniqueness
for parabolic differential equations. Firstly we establish local 
existence (with
respect to time) and then show that the solution remains bounded 
for any $T > 0$.
For (\ref{stkeq}) these steps are preceeded by the introduction of a change of
variable which enables us to consider, instead of the stochastic differential
equation, a related deterministic equation. Local existence and
 uniqueness is then
established via an application of a fixed point argument over a 
suitably defined 
space.

The application of the fixed point theorem necessitaties expressing 
the nonlinear 
solution operator of the derived deterministic equation as 
a mapping from a space
$\calE$ into itself. To achieve this we must show that the extensions of two
operators, which arise in the analysis, are well defined. This effort
is the major part of section \ref{sec-lclex}.

In section \ref{sec-glbex} we show that the local solution 
remains bounded for
any $T > 0$ which implies global existence of the solution.

We begin our discussion by presenting in the next section several 
definitions and 
basic results which we use later in our analysis.

\section{Preliminaries}
\label{sec-prel}

As usual, we denote by $L^{p}(I)$, $p = 1, 2, \ldots $ the closure of
$C^{\infty}(I)$ (the space of infinitely differentiable functions on $I$)
with respect to the $L^{p}(I)$ norm:
\[
   \| f \|_{L^{p}(I)} \ = \ \left( \int_{I} \ | f |^{p} \ dx \right)^{1/p} \ .
\]
Also, $H_{0}^{k}(I)$, $k = 1, 2, \ldots $, denotes the closure of 
$C_{0}^{\infty}(I)$ (the space of infinitely differentiable functions on $I$
which vanish at the the endponts)
with respect to the $H_{0}^{k}(I)$ norm:
\[
    \| f \|_{H^{k}_0(I)} \ = \ \left( \int_{I} \ | f |^{2} \ dx \ + \ 
      \int_{I} \ | f' |^{2} \ dx \ + \  \ldots \ + \ 
       \int_{I} \ | f^{(k)} |^{2} \ dx \  \right)^{1/2} \ .
\]
For convenience we use 
\[
        H \ := \ L^{2}(I) \ \mbox{ and } \ V \ := \ H_{0}^{1}(I) \, .
\]
 
To account for the temperal dependence we use the Banach spaces 
$L^{p}( 0 , T ; L^{q}( I ) )$, with the associated norm:
\[
    \| f \|_{L^{p}( 0 , T ; L^{q}( I ) )} \ := \ 
  \left( \int_{0}^{T} \, \left( \int_{I} \, | f |^{q} \, dx \right)^{p/q} 
           dt \right)^{1/p} \ .
\]

\textbf{Note}:
 The spaces $L^{p}( 0 , T ; H_{0}^{k}( I ) )$ are defined analogously.

A central role in the analysis below is played by the space $E$ defined by
\[
      E \ := \ L^{4}( 0 , T ; L^{4}( I ) ) \ .
\]
We remark that this choice for $E$ arises from the proof of lemma \ref{lmaGx} and
is dictated by the nonlinear term $u u_{x}$.

We begin by establishing the following embedding result which we combine with 
lemma \ref{lmaregy} to establish the setting for the application of the 
Banach contraction mapping theorem (lemma \ref{lmaBan}).
\begin{lemma}
 \label{lmainc}
 For any $T > 0$ we have
\begin{equation}
 L^{\infty}(0 , T; H) \cap L^{2}(0 , T; V) \subset E \ ,
 \label{eqinc}
\end{equation}
and there exists a constant $K$, independent of $T > 0$, such that
\begin{equation}
 \| u \|_{E} \le K \left( \| u \|_{L^{\infty}(0 , T; H)} \ + \ 
       \| u \|_{L^{2}(0 , T; V)}  \right) \ , \ u \in E \ .
 \label{eqinc2}
\end{equation}
\end{lemma}
\textbf{Proof}:
  We have by the Sobolev embedding theorem, (see \cite{ada751}, pg. 217), that
\[
     H^{1/2}(D) \subset  H^{1/4}(D) \subset L^{4}(D)
\]
and there exists a constant $C_{1} > 0$ such that 
\begin{equation}
 \| v \|_{L^{4}(D)} \  \le \ C_{1} \| v \|_{H^{1/2}(D)} \ .
\label{eqbd1}
\end{equation}

Using the interpolation inequality for ${H^{1/2}(D)}$ in terms
of $L^{2}(D)$ and ${H^{1}(D)}$ we have for
some constant $C_{2} > 0$ and all $t \in [0 , T]$
\begin{equation}
   \| u( t ) \|_{H^{1/2}(D)} \le C_{2} \, \| u( t ) \|_{L^{2}(D)}^{1/2} \ 
                     \| u( t ) \|_{H^{1}(D)}^{1/2} \ .
\label{eqbd2}
\end{equation}
Raising both sides of (\ref{eqbd2}) to the fourth power and
integrating (\ref{eqbd2}) over the interval $[0 , T]$, equation (\ref{eqinc2})
follows using standard inequalitites and the definitions of the norms,
with $K \ = \ C_{1} \, C_{2} / 2$ . \\
$\mbox{  }$ \hfill \qed

Essential to establishing the local existence is the following 
contraction mapping theorem.

\begin{lemma} (\cite{dap961}, Pg. 290)
 \label{lmaBan}
Let $\calF$ denote a transformation from a Banach space $\calE$ into
$\calE$, $\tilde{a}$ an element of $\calE$ and $\alpha > 0$ a positive number.
If $\calF(0) = 0$, $\| \tilde{a} \| \le \frac{1}{2} \alpha$ and
\begin{equation}
 \| \calF(z_{1}) - \calF(z_{2}) \|_{\mathcal{E}} \le \frac{1}{2} 
                                          \| z_{1} - z_{2} \|_{\mathcal{E}} 
      \mbox{ for } \| z_{1} \|_{\mathcal{E}} \le \alpha , \
                   \| z_{2} \|_{\mathcal{E}} \le \alpha ,
 \label{eqaF}
\end{equation}
then the equation
\begin{equation}
     z \ = \ \tilde{a} \ + \ \calF(z) , \ \ z \in \calE ,
 \label{eqfp}
\end{equation}
has a unique solution $z \in \calE$ satisfying $\| z \|_{\mathcal{E}} \le \alpha $.
\end{lemma}
$\mbox{  }$ \hfill \qed

Below we use the following lemma and corollary, 
which describe the regularity of the solution to a negative
self--adjoint operator.  
\begin{lemma} (\cite{zei901} pg. 424)
\label{lmaregy}
 Assume that $A$ is a negative self--adjoint operator on $H$ and 
\[
    V \ = \ D( ( -A )^{1/4} ) \subset H \subset V' \ .
\]
Then $A$ and $S(t) \ = \ e^{t A}$ has a continuous extension from 
$V$ to $V'$. If
\[
   y(t) \ = \ y(t ; g) \ = \ S( t ) y_{0} \ + \ 
                \int_{0}^{t} \, e^{(t - s) A} \, g(s) \, ds \ ,
             t \in [0 , T] , 
\]
for $y_{0} \in H$, and $ g \in L^{2}(0 , T; V') $,
then
\[
    y \in  L^{\infty}(0 , T; H) \cap L^{2}(0 , T; V) ,
\]
and for some constant $L > 0$, independent of $T > 0$,
\begin{equation}
    \| y \|_{L^{\infty}(0 , T; H)} \ + \ \| y \|_{L^{2}(0 , T; V)} \ 
  \le  \ L \left( \| y_{0} \|_{H} \ + \ \| g \|_{L^{2}(0 , T; V')} \right) \ .
\label{ctyS}
\end{equation}
\end{lemma}
$\mbox{  }$ \hfill \qed

\begin{corollary}
\label{correg}
For $A$, $S(t)$, and $y_{0}$ as described in lemma \ref{lmaregy}, we have
that
\begin{equation}
  \| S(t) y_{0} \|_{E} \le 8 K T^{1/4} 
  \left( \| S(t) y_{0} \|_{L^{\infty}(0, T; H)}^{4} \ + \
      \| S(t) y_{0} \|_{L^{\infty}(0, T; V)}^{4} \right)^{1/4} \ .
\label{estSy}
\end{equation}
\end{corollary}
\textbf{Proof}:  
Using (\ref{eqbd1}), (\ref{eqbd2}), (\ref{ctyS}) and, for notation convenience,
$u \ := \ S(t) y_{0}$, we have that
\begin{eqnarray*}
\int_{0}^{T} \ \| u \|_{L^4}^{4} \ dt &\le&
  ( C_{1} \, C_{2} )^{4} \ \int_{0}^{T} \ 
      \| u \|_{H}^{2} \ \| u \|_{V}^{2} \ dt  \\
 &\le& ( C_{1} \, C_{2} )^{4} / 2 \left(
     \int_{0}^{T} \  \| u \|_{H}^{4} \ dt  \ + \
    \int_{0}^{T} \  \| u \|_{V}^{4} \ dt  \right) \\
 &\le& ( C_{1} \, C_{2} )^{4} T / 2 \left(
      \  \| u \|_{L^{\infty}(0, T; H)}^{4}  \ + \
      \| u \|_{L^{\infty}(0, T; V)}^{4}  \right) \ .
\end{eqnarray*}
Taking the fourth root of both sides yields (\ref{estSy}) for
$K \ = \ C_{1} \, C_{2} / 2 $. \\
$\mbox{  }$ \hfill \qed

\textbf{Note}: From lemma \ref{lmaregy} we have that 
$S(t) \, y_{0} \in L^{\infty}(0, T; H) \cap L^{2}(0, T; V)$ which 
guarantees that the right hand side of (\ref{estSy}) is finite.

%
\section{Local Existence and Uniqueness}
\label{sec-lclex}
Our first step in establishing local existence and uniqueness of the stochastic
differential equation is to introduce a change of variable to reduce (\ref{stkeq})
to a deterministic equation.

Denote by $A$ the self-adjoint operator
\begin{equation}
  A u \ := \ - u_{xxxx} \ - u_{xx} \ - c \, u \ .
\label{defA}
\end{equation}
We assume that $c$ is chosen sufficiently large such that A is a strictly 
negative operator on the space $H^{4}_{0}( I )$.

Observe that as A is a strictly negative, self--adjoint, operator we can define
$(-A)^{\alpha}$ via Fourier analysis, with domain $D( (-A)^{\alpha} ) \ = \ 
H^{4 \alpha}_{0}( I )$. (See \cite{tem881} pg.55 for details.)

In view of (\ref{defA}) note that (\ref{stkeq}) can be rewritten in the form
\begin{equation}
  du \ = \ ( A u \ - u u_{x} \ + c u ) \ dt \ + \ dw \ ,
\label{stkmeq}
\end{equation}
where the Wiener process $w$ takes value in the separable Hilbert space 
$H = L^2(I)$ and it has the covariance operator $Q$.
With $S( t ) \ := \ e^{t A}$, $t \ge 0$, we define $w_{A}( t )$ via the 
stochastic integral
\begin{equation}
      w_{A}( t ) \ := \ \int_{0}^{t} S( t - s ) \, dw(s) \ .
\label{defwa}
\end{equation}
Using the substitution
\begin{equation}
   y(t , x) \ := \ u(t , x) \ - \ w_{A}(t , x) \ , \ t \in [0 , T] \ \prob \mbox{--} a.s. \ ,
\label{defy}
\end{equation}
(\ref{stkeq}) reduces to the deterministic problem
\begin{equation}
 y_{t} \ = \ A y - (y \ + w_{A}) (y \ + w_{A})_{x} \ + c (y \ + w_{A}) \ ,
\label{deteq}
\end{equation}
subject to 
\begin{equation}
 y(0 , x) \ = \ u_{0}(x) \ \mbox{ and } y(t , -l) \ = \ y(t , l) \ = \ 0 \ .
\label{bcy}
\end{equation}

\textbf{Note}: 
The assumption that $dw$ in (\ref{stkeq}) denotes a Q--Wiener process,
together with the fact that $A$ is a strictly negative self--adjoint operator,
ensures that $w_{A}(t)$ given by (\ref{defwa}) 
has a version which is H\"{o}lder continuous
with values in $D( (-A)^{\alpha} )$ for $0 \le \alpha < 1/4$, with 
H\"{o}lder exponent less than $(1/4 - \alpha)$, (see \cite{dap961}, pg. 60).
Thus, with $\alpha \ = \ 1 / 8$, in view of (\ref{eqbd1}), we conclude
that $w_{A}(t)$ has a continuous version in $L^{4}(I)$.
Below we take $w_{A}(t)$ to denote this continuous version.

The solution $y$ satisfying (\ref{deteq}) may be expressed in integral form
as
\begin{eqnarray}
 y(t) & = & S( t ) u_{0} \ + \ \nonumber \\
   & &  \hspace{-0.3in} \int_{0}^{t} S( t - s ) \left[ 
   - ( y( s ) \ + \ w_{A}( s ) ) ( y( s ) \ + \ w_{A}( s ) )_{x} \ + 
    \ c ( y( s ) \ + \ w_{A}( s ) ) \right] \ ds   
     \label{eqint1} \\
      & = & S( t ) u_{0} \ + \ F( y \ + \ w_{A} )( t ) \ , \ \ t \in [0 , T] \ .
        \label{eqint2}
\end{eqnarray}
In the following we show the existence and uniqueness of the solution
$y$ to this integral equation (\ref{eqint2}). This gives a so-called
(mild) solution $u$ for the stochastic
Kuramoto-Sivashinsky equation (\ref{stkeq}). This is the definition of  
`solution' used in this paper.

In (\ref{eqint2}) $F : E \rightarrow E $ is a continuous extension of the operator 
\[
   F_{0} : C^{1}([0 , T] ; V) \rightarrow E
\]
defined by
\begin{equation}
  ( F_{0} u )( t ) \ = \ \int_{0}^{t} S(t - s) ( G_{0} u )( s ) \ ds \ , 
    t \in [0 , T] ,
 \label{defF0}
\end{equation}
where
\[
   G_{0} : C^{1}([0 , T] ; V) \rightarrow E
\]
is given by
\begin{equation}
  ( G_{0} u )( t ) \ = \ - u u_{x}(t) \ + \ c u(t) \ , 
    t \in [0 , T] .
 \label{defG0}
\end{equation}
  
In view of (\ref{eqint2}), and assuming that $F$ is well defined --- a non--trivial point 
whose discussion occupies the later part of this section --- , on applying lemma
\ref{lmaBan}
we have local existence of the solution to (\ref{stkeq}),(\ref{bc1}).

\textbf{Note}:
 The value of $\tau$ for which we establish local existence and uniqueness of the
solution on the interval $ [0 , \tau] $, depends upon the particular realization.

\begin{theorem} 
 \label{locex}
 For $u_{0}$ in $H$ there exists a random variable $\tau$ taking values $\prob$--$a.s.$
in $(0 , T]$ such that equations (\ref{stkeq})(\ref{bc1}) 
have a unique solution $u$ on
the interval $[0 , \tau]$.
\end{theorem}

Note that by a general result in \cite{dap921}, page 72, the solution 
$u$ has a measurable modification. In the following the solution $u$ refers
to this measurable version.

\textbf{Proof}:
Observe that with $z(t) \ = \ y(t) \ + \ w_{A}(t) \ - \ S(t) \, u_{0} $,
equation (\ref{eqint2}) may be rewritten as
\begin{equation}
     z \ = \ \tilde{a} \ + \ \calF(z) 
\label{defzeq}
\end{equation}
for $ \tilde{a} \ = \ w_{A}(t)$, and 
$\calF(z) \ = \ F(z \ + \ S(t) \, u_{0})$.
Thus the existence and uniqueness of the solution to (\ref{eqint2})
is equivalent to that for (\ref{defzeq}). 

Let $\alpha = 1/6 M$, and $\tau_{1}$ be
given by
\begin{equation}
  \tau_{1} \ = \ \left(6 M \, [ c \, (2 l)^{1/4} \ + \ 
  16 K \, ( \| S(t) u_{0} \|_{L^{\infty}(0, T; H)}^{4} \ + \ 
     \| S(t) u_{0} \|_{L^{\infty}(0, T; V)}^{4}  )^{1/4} \, ] \right)^{-4} \ ,
 \label{dt1}
\end{equation}
for $K$ defined in lemma \ref{lmainc}, and $M$ defined in lemma \ref{lmextF}.

The $\tau_{1}$ is well-defined and it is so chosen  that it will

guarantee that $\cal F$ is a contraction mapping (see below).

As $w_{A}(t)$ is continuous with $w_{A}(0) \ = \ 0$, there exists
$\tau_{2}$ such that
\[
   \int_{0}^{t} \, \| w_{A}(s) \|_{L^{4}(I)}^{4} \, ds  \ \le \ \alpha / 2 \, , 
 \ \mbox{ for } 0 \le t \le \tau_{2} \ .
\]

Let $\tau := \min \{\tau_{1} \, , \, \tau_{2} \}$ and analogous to the
definition for $E$ introduce $\calE$ as
\[
        \calE \ := \ L^{4}(0 , \tau ; L^{4}(I)) \ .
\]

With $z_{1}$ and $z_{2}$ satisfying $\| z_{i} \|_{\calE} \le \alpha $ ( = 
  1 / 6 M ) for $i = 1 , 2$, we have using lemma \ref{lmextF}, 
(\ref{estSy}), and the definition of $\tau$
\begin{eqnarray*}
 \| \calF(z_{1}) - \calF(z_{2}) \|_{\calE} &\le&
        M \left( \| z_{1} \ + \ S(t) u_{0} \|_{\calE} \ + \
           \| z_{2} \ + \ S(t) u_{0} \|_{\calE} \ + \ 
              c ( 2 l \tau )^{1/4} \right) \ 
           \| z_{1} \ - \ z_{2} \|_{\calE}  \\
   &\le& M \left( \| z_{1} \|_{\calE} \ + \
           \| z_{2} \|_{\calE} \ + \ 2 \| S(t) u_{0} \|_{\calE} \ + \ 
              c ( 2 l \tau )^{1/4} \right) \ 
           \| z_{1} \ - \ z_{2} \|_{\calE}  \\
    &\le& \frac{1}{2} \| z_{1} \ - \ z_{2} \|_{\calE} \ .
\end{eqnarray*}

Finally, applying lemma \ref{lmaBan} we establish the existence and
uniqueness of $z(t)$, and consequently $y(t)$, on the 
interval $[0 , \tau]$. \\
$\mbox{  }$ \hfill $\qed$

What remains is to establish the regularity result used for $F$ in the proof
of theorem \ref{locex}. However we must first show that $F$ is well defined by 
showing $G_{0}$ and $F_{0}$ defined by (\ref{defG0}) and (\ref{defF0}) have 
appropriate extensions.

\begin{lemma}
 \label{lmaGx}
The operator $G_{0}$ defined by (\ref{defG0}) can be continuously extended to
\[
    G : \ E \ \rightarrow \ L^{2}(0 , T ; V') \ ,
\]
satisfying
\[
    \| G(u) \ - \ G(v) \|_{L^{2}(0 , T ; V')} \le
               27^{1/4} ( \| u \|_{E} \ + \| v \|_{E} \ + c (2 l T)^{1/4} ) 
                    \| u - v \|_{E}  \ \ u , v \in E \ .
\]
\end{lemma}
\textbf{Proof}:
Let $u, \ v, \ \psi \in L^{2}(0 , T; V)$. Denoting the 
duality mapping between $L^{2}(0 , T; V)$ and $L^{2}(0 , T; V')$ by
$\langle \cdot , \cdot \rangle$, we have
\begin{eqnarray*}
 \langle G_{0}( u ) - G_{0}( v ) , \psi \rangle & = & \int_{0}^{T} \int_{I} \left( 
       - u u_{x} \ + v v_{x} \right) \psi \ + \ c (u - v) \psi \ dx \ dt \\
         & = & \int_{0}^{T} \int_{I} \frac{1}{2} (u + v)(u - v) \psi_{x} \ + \
           c (u - v) \psi \ dx \ dt \ ,
\end{eqnarray*}
i.e.
\begin{eqnarray*}
| \langle G_{0}( u ) - G_{0}( v ) , \psi \rangle | & \le & 
    \left( \int_{0}^{T} \int_{I} (|u| \ + \ |v| + \ c)^{2}(u \ - \ v)^{2} 
     \ dx \ dt \right)^{1/2} \ \cdot \\
    & & \hspace{1.0in} 
    \left( \int_{0}^{T} \int_{I} (\frac{1}{2} |\psi_{x}| \ + \ |\psi|)^{2} 
     \ dx \ dt \right)^{1/2} \\
  & \le &
    \left( \int_{0}^{T} \int_{I} (|u| \ + \ |v| + \ |c|)^{4} \ dx \ dt \right)^{1/4} \ 
      \cdot \\
    & & \hspace{0.5in} 
    \left( \int_{0}^{T} \int_{I} (u \ - \ v )^{4} \ dx \ dt \right)^{1/4} \ 
    \| \psi \|_{L^{2}(0 , T ; V)}    \\
  & \le &
    \left( 27 ( \| u \|_{E}^{4} \ + \ \| v \|_{E}^{4} \ + \ \| c \|_{E}^{4} )
      \right)^{1/4} \cdot \\
    & & \hspace{1.0in} 
      \|(u \ - \ v )\|_{E} \ \| \psi \|_{L^{2}(0 , T ; V)}    \\
  & \le &
     27^{1/4} ( \| u \|_{E} \ + \ \| v \|_{E} \ + \  c (2 l T)^{1/4} )  \cdot \\
    & & \hspace{1.0in} 
       \ \|(u \ - \ v )\|_{E} \
     \| \psi \|_{L^{2}(0 , T ; V)}   \ ,
\end{eqnarray*}
from which the result follows.
$\mbox{  }$ \hfill $\qed$

For the extension of $F_{0}$ and the regularity of $F$ we have:
\begin{lemma}
\label{lmextF}
The transformation $F_{0}$ described by (\ref{defF0}) can be continuously extended
to $F : \ E \rightarrow E$. Moreover there exists a constant $M > 0$, independent of
$T > 0$, such that
\begin{equation}
  \| F(u) \ - \ F(v) \|_{E} \le M  ( \| u \|_{E} \ + \ \| v \|_{E} \ + 
      \  c (2 l T)^{1/4} ) \|(u \ - \ v )\|_{E} \ , \ u, \, v \in E \ .
\label{regF}
\end{equation}
\end{lemma}

\textbf{Proof}:
In view of the definition of $F_{0}$ in (\ref{defF0}) and lemma \ref{lmaregy} we have that
\[
    F( u ) (t) \ = \ y(t ; G(u)) \in L^{\infty}(0 , T; H) \cap L^{2}(0 , T; V) \ .
\]
Moreover from lemma \ref{lmainc} we have
\begin{eqnarray*}
 \| F(u) \ - \ F(v) \|_{E} & \le &
             K \left( \| y(\cdot ; G(u)) \ - \ y(\cdot ; G(v)) \|_{L^{\infty}(0 , T; H)}
        \ +  \right. \\
 & & \hspace*{1.5in} \left. 
          \| y(\cdot ; G(u)) \ - \ y(\cdot ; G(v)) \|_{L^{2}(0 , T; V)} \right) \\
     & \le & 
        K L \| G(u) - G(v) \|_{L^{2}(0 , T; V')}  \ , \mbox{ using lemma \ref{lmaregy} } , \\
    & \le & M  ( \| u \|_{E} \ + \ \| v \|_{E} \ + 
      \  c (2 l T)^{1/4} ) \|(u \ - \ v )\|_{E} \ , \ \ \mbox{ using lemma \ref{lmaGx} } ,
\end{eqnarray*}
for $M  = 27^{1/4} K L $ . \\
$ \mbox{ } $ \hfill $\qed$

%
\section{Global Existence}
\label{sec-glbex}
We now extend the local existence of theorem \ref{locex} established in the previous
section to global existence. Local existence establishes that the solution, $u$,
lies in the solution space, $E$, for some initial time period. We establish 
global existence by showing that for any time $T$ the $E$--norm of $u$ is finite
and hence $u$ still lies in the space $E$.
To do this we first establish that the solutions are continuous with respect
to the initial data. This enables us to restrict our attention to showing
that \underline{strong solutions} remain bounded in $E$.

Following directly the proof of lemma \ref{lmaGx} with the inner product 
(only) taken over the spatial domain, $I$, we have:

\begin{corollary} \label{corGnmVp}
The operator $G$ defined in lemma \ref{lmaGx} satisfies for 
$\ u , \ v \in E$ 
\begin{equation}  \label{GnmVp}
 \| G(u) \ - \ G(v) \|_{V'} \le  27^{1/4} 
   ( \| u \|_{L^{4}(I)} \ + \| v \|_{L^{4}(I)} \ + c (2 l)^{1/4} ) \
                    \| u - v \|_{_{L^{4}(I)}}  \ .
\end{equation}
\end{corollary}
 \mbox{  } \hfill \qed

\begin{lemma} \label{ctdep}
The solution, $y(t)$, of (\ref{eqint2}) depends continuous on the initial
data $u_{0} \in H$, and the random forcing term $w_{A}(t) \in E$.
\end{lemma}

Note that the continuous dependence on $w_{A}(t)$ is needed in the proof of the
next lemma where we approximate $w_{A}(t)$ by regular processes.

\textbf{Proof}:
Let $y_{0}$, and $y_{1}$ denote solutions of (\ref{eqint2}) generated by
 $u_{0}$, $w_{A}^{0}(t)$, and $u_{1}$, $w_{A}^{1}(t)$, respectively. Then,
 on the common existence interval of $y_{0}$ and $y_{1}$,
\begin{eqnarray*}
 y_{0} \ - \ y_{1} & = & S( t ) ( u_{0} \ - \ u_{1} ) \ + \  \\
   & &   \int_{0}^{t} \left[  S( t - s ) G(y_{0} + w_{A}^{0})(s)
   \ - \ S( t - s ) G(y_{1} + w_{A}^{1})(s) \right] \ ds   \ .
\end{eqnarray*}
From lemma \ref{lmaregy} we have $y_{0}$, and $y_{1}$ $\in L^{2}(0,T;V)$
thus $(y_{0} \ - \ y_{1} )(t) \in V , \mu \ a.e.$, (i.e. for almost all
$t \in (0 , T)$). Using the continuity of $S(t)$ and lemma \ref{corGnmVp} we
have that there exits constants $L_{1}$ and $C_{1}$ such that
\begin{eqnarray*}
 \| y_{0} \ - \ y_{1} \|_{V} & \le & L_{1} \| u_{0} \ - \ u_{1} \|_{H} \ + \
    \int_{0}^{t} C_{1} \| G(y_{0} + w_{A}^{0})(s)
   \ - \  G(y_{1} + w_{A}^{1}) \|_{V'}  \ ds   \\
  & \le & L_{1} \| u_{0} \ - \ u_{1} \|_{H} \ + \ \\
  & & \hspace*{-1.0in} C_{1} \int_{0}^{t}  27^{1/4} 
   ( \| y_{0} + w_{A}^{0} \|_{L^{4}(I)} \ + \| y_{1} + w_{A}^{1} \|_{L^{4}(I)} \ +
   c (2 l)^{1/4} ) 
        \| ( y_{0} + w_{A}^{0} ) \ - \ (  y_{1} + w_{A}^{1} ) \|_{_{L^{4}(I)}} \ ds \ , \\
 & & \mu \ a.e.
\end{eqnarray*}
Using the Sobolev embedding theorem, the existence of $y_{0}$, and 
$y_{1}$ $\in E$, implies there exists constants $C_{2}$ and $C_{3}$
such that
\[
 \| y_{0} \ - \ y_{1} \|_{L^{4}(I)}   \le   
 C_{2} \left( \| u_{0} \ - \ u_{1} \|_{H} \ + 
                  \ \| w_{A}^{0} \ - \ w_{A}^{1} \|_{E} \right) \ + \ 
    C_{3} \int_{0}^{t} \| y_{0} \ - \ y_{1} \|_{L^{4}(I)} \ ds \ \ \mu \ a.e.
\]
Applying Gronwall's inequality then yields
\[
\| y_{0} \ - \ y_{1} \|_{L^{4}(I)}   \le  
    C_{2} \left( \| u_{0} \ - \ u_{1} \|_{H} \ + 
                  \ \| w_{A}^{0} \ - \ w_{A}^{1} \|_{E} \right)
 \ e^{C_{3} t} \ \ \ \mu \ a.e.
\]
from which the stated conclusion follows.
 \mbox{  } \hfill \qed

Next we establish appropriate norm estimates for the solution.

\begin{lemma}
\label{glnmest}
Let $u_{0} \in H$, $w_{A}$ be given by (\ref{defwa}), and $y$ denote the 
solution of 
\begin{equation}
   y(t) \ = \ S( t ) u_{0} \ + \ F( y \ + \ w_{A} )( t ) \ , \ \ t \in [0 , T] \ .
\label{eqgex1}
\end{equation}
Then, $y$ satisfies
\begin{equation}
  \sup_{t \in [0 , T]} \, \| y( t ) \|_{H}^{2} \le 
   \| u_{0} \|^2_{H} \, e^{\int_{0}^{T} \, f( s ) \, d s} \ + \ 
   \int_{0}^{T} \, e^{\int_{0}^{T} \, f( s ) \, d s} \ g( \tau ) \, d \tau \ ,
 \label{Hnmest}
\end{equation}
and
\begin{equation}
 \int_{0}^{T} \, \| y ( t ) \|_{V}^{2} \, dt  \le  \| u_{0} \|^2_{H} \ + \
    \sup_{s \in [0 , T]} \, \| y( s ) \|_{H}^{2} \, 
    \int_{0}^{T} \, f( \tau ) \, d \tau \
    + \ \int_{0}^{T} \, g( \tau ) \, d \tau  \ ,
 \label{Vnmest}
\end{equation}
where 
\begin{eqnarray*}
  f( t ) & = & \frac{1}{2} \left( 2 \, c \ + \ (2 \, + \, \frac{3}{4} C_{1} C_{2})^{2} \ + \ 
      C_{1} C_{2} \| w_{A} \|_{L^{4}}^{4} \right) \ , \\
 \mbox{ and } & &  \hfill  \\
  g( t ) & = & c \, \| w_{A} \|_{H}^{2} \ + \ 
            \frac{1}{4} \| w_{A} \|_{L^{4}}^{4} \ .
\end{eqnarray*}
\end{lemma}

\textbf{Proof}:
As $D( A )$ and $C(0, T; H_{0}^{2}(I))$ are dense in $H$
and $E$, respectively, and from lemma \ref{ctdep}
we have established continuous dependence of the solution, 
it suffices to establish (\ref{Hnmest}),(\ref{Vnmest}) 
for the (strong) solution of the differential equation
\begin{eqnarray}
 \frac{d y(t)}{d t} & = & A y(t) \ - \ (y(t) \ + \ w_{A}(t)) (y(t) \ + \ w_{A}(t))_{x} \
 + \ c (y(t) \ + \ w_{A}(t))  \ ,  \label{eqde1} \\
 y_{0} & = & u_{0} \ . \nonumber 
\end{eqnarray}

We first show that 
\begin{equation}
 \frac{d}{d t} \| y(t) \|_{H}^{2} \ + \ \| y(t) \|_{V}^{2} \le f(t) \, \| y \|_{H}^{2} \
 + \ c \, \| w_{A} \|_{H}^{2} \ + \ \frac{1}{4} \| w_{A} \|_{L^{4}}^{4} \ .
\label{nmieq1}
\end{equation}

Multiplying (\ref{eqde1}) by $y(t)$ and integrating over $I$ we obtain
\begin{equation}
 \frac{1}{2} \frac{d}{d t} \| y \|_{H}^{2} \ + \ \| y_{xx} \|_{H}^{2}  \ - \ 
    \| y_{x} \|_{H}^{2} \ = \ 
 - \int_{-l}^{l} \, y (y \ + \ w_{A}) (y \ + \ w_{A})_{x} \, dx \ + \ 
   c \int_{-l}^{l} \, y w_{A} dx \ .
 \label{nmieq2}
\end{equation}

We obtain a lower bound for $\| y_{xx} \|_{H}$ via:
\begin{eqnarray*}
\| y_{x} \|_{H}^{2} & = & \int_{-l}^{l} y_{x} \, y_{x} \, dx \ = \ 
  -\int_{-l}^{l} \, y_{x x} \, y \, dx  \\
     & \le & \left( \int_{-l}^{l} \, y_{x x}^{2} \, dx \right)^{1/2} \,
         \left( \int_{-l}^{l} \, y^{2} \, dx \right)^{1/2}  \\
     & \le & \frac{1}{(2 \, + \, \frac{3}{4} C_{1} C_{2})} \, \| y_{xx} \|_{H}^{2} \
 + \ \frac{(2 \, + \, \frac{3}{4} C_{1} C_{2})}{4} \, \| y \|_{H}^{2} \ .
\end{eqnarray*}
Rearranging yields
\begin{equation}
 \| y_{xx} \|_{H}^{2} \ \ge \ (2 \, + \, \frac{3}{4} C_{1} C_{2}) \, \| y_{x} \|_{H}^{2} \
  - \ \frac{(2 \, + \, \frac{3}{4} C_{1} C_{2})^{2}}{4} \, \| y \|_{H}^{2} \ .
 \label{nmieq3}
\end{equation}

Next consider the first integral on the r.h.s. of (\ref{nmieq2}) in three pieces:
\begin{eqnarray}
 -\int_{-l}^{l} \, y \, y \, y_{x} \, dx & = & -\int_{-l}^{l} \, \frac{1}{3} 
    \frac{d}{d x}( y^{3} ) \, dx \ = \ 0 \ , 
  \label{nmieq4}  \\
 | -\int_{-l}^{l} \, y \, w_{A} \, {w_{A}}_{x} \, dx \, | & = & 
 | \frac{1}{2} \int_{-l}^{l} \, y_{x} \, w_{A}^{2} \, dx \, | \ \le 
 \frac{1}{2} \| y \|_{V}^{2} \ + \ \frac{1}{8} \| w_{A} \|_{L^{4}}^{4} \ ,
   \label{nmieq5}  \\
 \mbox{ and } \hspace*{2in} & &  \nonumber \\
  -\int_{-l}^{l} \, \left(y \, y_{x} \, w_{A} \ + \ y^{2} \, {w_{A}}_{x} \right) \, dx \,
    & = &
  \int_{-l}^{l} \, y \, y_{x} \, w_{A} \, dx \,  .
  \nonumber
\end{eqnarray}
This last term is estimated as follows.
\begin{eqnarray}
 | \int_{-l}^{l} \, y \, y_{x} \, w_{A} \, dx | & \le &
       \left( \int_{-l}^{l} \, y_{x}^{2} \, dx \right)^{1/2} \, 
        \left( \int_{-l}^{l} \, y^{2} \, w_{A}^{2} \, dx \right)^{1/2} \, 
   \nonumber \\
   & \le &
       \left( \int_{-l}^{l} \, y_{x}^{2} \, dx \right)^{1/2} \, 
        \left( \int_{-l}^{l} \, y^{4} \, dx \right)^{1/4} \, 
         \left( \int_{-l}^{l} \, w_{A}^{4} \, dx \right)^{1/4} \, 
    \nonumber  \\
   & \le &
          \| y \|_{V} \, \| y \|_{L^{4}} \, \| w_{A} \|_{L^{4}}
   \nonumber \\
  & \le &
          C_{1} \, C_{2} \| y \|_{V}^{3/2} \, 
           \| y \|_{H}^{1/2} \,\| w_{A} \|_{L^{4}} \mbox{ ( using (\ref{eqbd1}),(\ref{eqbd2}) )} 
   \nonumber \\
  & \le &
          \frac{3 \, C_{1} \, C_{2}}{4} \| y \|_{V}^{2} \ + \
           \frac{C_{1} \, C_{2}}{4}\| y \|_{H}^{2} \,\| w_{A} \|_{L^{4}}^{4}
                 \label{nmieq7} \\
   & & \hspace*{2in}       \mbox{ ( using Young's Inequality )} .
         \nonumber
\end{eqnarray}

For the remaining term on the r.h.s. of (\ref{nmieq2}) we have
\begin{eqnarray}
 | c \int_{-l}^{l} \, y \, w_{A} \, dx | & \le &
      c \left( \int_{-l}^{l} \,  y^{2} \, dx \right)^{1/2} \,
        \left( \int_{-l}^{l} \,  w_{A}^{2} \, dx \right)^{1/2} 
     \nonumber \\
     & \le &
      \frac{c}{2} \| y \|_{H}^{2} \ + \  \, \frac{c}{2} \| w_{A} \|_{H}^{2} \ .
    \label{nmieq8}
\end{eqnarray}

Combining (\ref{nmieq2})---(\ref{nmieq8}) yields (\ref{nmieq1}).

The estimate for $\| y \|_{H}$ in (\ref{Hnmest}) now follows from the observation
that (\ref{nmieq1}) is a first order differential inequality for $\| y \|_{H}$. 
The bound involving $\| y \|_{V}$ in (\ref{Vnmest}) is established by integrating 
(\ref{nmieq1}) from
$0$ to $T$. \\
 \mbox{  } \hfill \qed

We are now in a position to establish the global existence of the solution.

\begin{theorem}
\label{glbex}
For $u_{0} \in H=L^2(I)$, there exists $\prob$ a.s. a unique solution 
$u(\cdot , x) \in E$ of (\ref{stkeq}),(\ref{bc1}). 
\end{theorem}

\textbf{Proof}:
From theorem \ref{locex} we have existence of the solution $u(\cdot , x)
\in E$, $\prob$ a.s., for the interval $[0 , \tau]$. In view of
(\ref{defy}) and lemma \ref{glnmest} we conclude that $u(t , x)$ remains
bounded in $E$, $\prob$ a.s. for all $t \ge 0$, which implies global
existence of the solution to (\ref{stkeq}),(\ref{bc1}). 
 \mbox{  } \hfill \qed

\bigskip

Finally we remark that by following the same argument as in
Brannan et al. \cite{Duan-SQG}, we can show that the solution
is actually H\"{o}lder continuous in space with exponent less
than $\frac18$. It is also possible to consider
multiplicative noise in equation (1.1). The approach in this paper should also
apply to other similar parabolic type stochastic
partial differential equations.




\end{document}